\documentclass[11pt]{article}
    
\usepackage{amsmath,amsfonts,amsthm,amscd,amssymb,graphicx}

\usepackage{subfigure}

\numberwithin{equation}{section}
\newtheorem{theorem}{Theorem}[section]

\newtheorem{remark}[theorem]{Remark}

\newcommand{\RR}{\mathbb{R}}

\newcommand{\TT}{{\mathbb T}}
\def\beq{\begin{equation}}
\def\eeq{\end{equation}}
\def\bb1{{1\!\!1}}
\def\rit{{\Bbb R}}

\def\tit{{\Bbb T}}

\def\app{\mathrm{app}}

\begin{document}

\title{Sublayer of Prandtl boundary layers}

\author{Emmanuel Grenier\footnotemark[1]
  \and Toan T. Nguyen\footnotemark[2]
}

\maketitle

\renewcommand{\thefootnote}{\fnsymbol{footnote}}

\footnotetext[1]{Equipe Projet Inria NUMED,
 INRIA Rh\^one Alpes, Unit\'e de Math\'ematiques Pures et Appliqu\'ees., 
 UMR 5669, CNRS et \'Ecole Normale Sup\'erieure de Lyon,
               46, all\'ee d'Italie, 69364 Lyon Cedex 07, France. Email: Emmanuel.Grenier@ens-lyon.fr}

\footnotetext[2]{Department of Mathematics, Penn State University, State College, PA 16803. Email: nguyen@math.psu.edu.} 
%TN's research was partly supported by the NSF under grant DMS-1405728.}

\begin{abstract} 

The aim of this paper is to investigate the stability of Prandtl boundary layers in the vanishing viscosity limit: $\nu \to 0$. 
In \cite{Grenier}, one of the authors
proved that there exists no asymptotic expansion involving one Prandtl's boundary layer with thickness of order $\sqrt\nu$, 
which describes the inviscid limit of Navier-Stokes equations. 
The instability gives rise to a viscous boundary sublayer whose thickness is of order $\nu^{3/4}$. 
In this paper, we point out how the stability of the classical Prandtl's layer is linked to the stability of this sublayer. 
In particular, we prove that the two layers cannot both be nonlinearly stable in $L^\infty$. 
That is, either the Prandtl's layer or the boundary sublayer is nonlinearly unstable in the sup norm. 

\end{abstract}

%%%%%%%%%%%%%%%%%%%%%%

\section{Introduction}

%%%%%%%%%%%%%%%%%%%%%%

In this paper, we are interested in the inviscid limit $\nu \to 0$ of the Navier-Stokes equations for incompressible fluids, namely
\beq \label{NS1}
\partial_t u^\nu + (u^\nu \cdot \nabla) u^\nu + \nabla p^\nu  = \nu \Delta u^\nu,
\eeq
\beq \label{NS2} 
\nabla \cdot u^\nu = 0,
\eeq
on the half plane $\Omega = \{(x,y) \in \tit \times \rit^+\}$ or the half space $ \Omega = \{(x,y) \in \tit^2 \times \rit^+\}$, with the no-slip boundary condition
\beq \label{NS3} 
u^\nu = 0 \quad \hbox{on} \quad \partial \Omega.
\eeq
%and initial condition
%\beq \label{NS4}
%u^\nu(0,\cdot) = u^\nu_0 .
%\eeq
As $\nu$ goes to $0$, one would expect the solutions $u^\nu$ to converge to solutions of Euler equations for incompressible fluids
\beq \label{Euler1}
\partial_t u^0 + (u \cdot \nabla) u^0 + \nabla p^0 = 0,
\eeq
\beq \label{Euler2}
\nabla \cdot u^0 = 0,
\eeq
with the boundary condition 
\beq \label{Euler3}
u^0 \cdot n = 0 \quad \hbox{on} \quad \partial \Omega,
\eeq
where $n$ is the unit normal to $\partial \Omega$.

At the beginning of the twentieth century, Prandtl introduced its well known boundary layers in order to describe the
transition from Navier Stokes to Euler equations as the viscosity tends to zero. Formally, we expect that
\beq \label{Ansatz}
u^\nu(t,x,y) \approx  u(t,x,y) + u_P \Bigl(t,x, {y \over \sqrt{\nu}}\Bigr) 
\eeq
where $u_P$ is the Prandtl boundary layer correction, which is of order one in term of small viscosity, and is having the boundary layer variable $y$ of order $\sqrt{\nu}$, the classical
size of Prandtl's boundary layer.

Prandtl boundary layers have been intensively studied in the mathematical literature. First, Oleinik \cite{Oleinik1,Oleinik} proved the existence
in small time of Sobolev solutions provided the initial vorticity is monotonic in the normal variable $z$. The Oleinik's monotonic solutions are also recently reconstructed via energy methods \cite{Xu, MW,VV1}. There are also analytic solutions to the Prandtl equations; see, for instance, \cite{Caf1, DGV1, VV} and the references therein. On the other hand, the authors in \cite{WE} construct a class of solutions which blow up in finite time. We also refer to \cite{VanDommelen,Gargano,VV2} for the study of the onset of singularities in Prandtl's
equations. Then, \cite{DGV} showed that Prandtl equations are ill posed in Sobolev
spaces for some classes of initial data; see also \cite{Grenier, GVN, GuoN,GGN1,GGN3} for further instability of Prandtl boundary layers.

The validity of Prandtl's Ansatz \eqref{Ansatz} has been established in \cite{Caf1,Caf2}
for initial data with analytic regularity, leaving a remainder of order $\sqrt \nu$. A similar result is also obtained in \cite{Mae}. If we assume only Sobolev regularity of the remainder in the approximation \eqref{Ansatz}, one of the authors proved in \cite{Grenier} that such an asymptotic expansion is false, up to a remainder of order $\nu^{1/4}$.

In this paper, we continue the analysis introduced in \cite{Grenier} to further study the structure of the instability of Prandtl's layers. Our aim is to analyze the boundary sublayer which prevents the previous analysis (\cite{Grenier}) to reach to instability of order one in its amplitude in the approximation \eqref{Ansatz}.

More precisely, we study the classical stability problem of whether the following time-dependent shear layer flow 
\begin{equation}\label{def-Pr}
U_s(t,x,y) = \Bigl( \begin{array}{c} U_s(t,y / \sqrt\nu) \\0 \end{array} \Bigr)
\end{equation}
is nonlinearly stable to the Navier-Stokes equation in the inviscid limit. Here, $U_s(t,z)$ solves the heat equation 
$$
\partial_t U_s = \partial_z^2 U_s
$$
with initial data $U_s(0,z) = U(z)$. Occasionally, we write $U$ in place of the vector $[U,0]^{tr}$. 

\bigskip

\noindent {\bf Assumption on $U_s$.}

\medskip

\noindent  We assume that the initial shear layer $U(z)$ is smooth, $U(0)=0$, and $\lim_{z\to \infty}U(z)$ is finite.
In addition, we assume that $U(z)$ is spectrally unstable to Euler equations. Precisely, there exists a growing solution of the form
\begin{equation}\label{def-vgr}
v_c(t,x,y) = v_s(x,y) e^{\Re \lambda t}
\end{equation}
solving the linearized Euler equations
$$
\partial_t v + (U \cdot \nabla) v + (v \cdot \nabla) U + \nabla p = 0,
$$
$$
\nabla \cdot v = 0 ,
$$
with the boundary condition $v \cdot n = 0$ on $\partial \Omega$. 

\medskip

Using the growing mode $v$, we will establish a nonlinear instability result for the classical $O(\nu^{1/2})$ Prandtl's layer $U_s$. 
This construction will involve a boundary sublayer of size $O(\nu^{3/4})$. To leading order (see Section \ref{sec-sublayer}), the sublayer is of the form 
\begin{equation}\label{def-sublayer1} v_S^1 = v_S^1 \Big(\frac{t}{\sqrt \nu},\frac{x}{\sqrt \nu},\frac{y}{\nu^{3/4}}\Big)\end{equation}
with $v_S^1$ solving the Stokes problem 
\begin{equation}\label{Stokes0}
\partial_t v_S^1 + \nabla p = \nu \Delta v_S^1 ,
\qquad \nabla \cdot v_S^1 = 0,
\end{equation}
with the following boundary conditions
$$ {v_S^1 }_{\vert_{z=0}}= 0, \qquad \lim_{z\to \infty} v_S^1(t,x,z) = v_c (t,x,0).$$
That is, the sublayer $v_S^1$ corrects the nonzero boundary condition of the inviscid growing mode $v_s$, defined in \eqref{def-vgr}. 
As it will be clear in the construction, 
$$
v_s = \Re \nabla^\perp (e^{i\alpha x} \psi_e(z))
$$
 with $\psi_e$ solving the Rayleigh equations with the zero Dirichlet boundary conditions.
 As a consequence, the boundary value of the tangential component of $v_s$ is nonzero, and the boundary sublayer is present in the construction. 
 
 %One observes that the sublayer $U_\sub$ has no inflection point, and hence $U_\sub(0,z)$ is spectrally stable to the corresponding Euler equations.  

Roughly speaking, we will prove that {\em the Prandtl's layer and the boundary sublayer cannot be simultaneously nonlinearly stable in $L^\infty$.} Precisely, we obtain the following theorem.

%\noindent The precise statement is technical and is postponed till the end of the paper.

\begin{theorem} Let $U_s$ be a Prandtl's boundary layer of the form \eqref{def-Pr}. Assume that the initial shear layer $U$ is spectrally unstable to the Euler equations, giving rise to the boundary sublayer $v_S^1$, defined as in \eqref{def-sublayer1}.   
Then, one of the following must hold 

 \begin{itemize}

\item For any $s,N$ arbitrarily large, we can find $\sigma_0 > 0$, initial
conditions $u^\nu(0)$ and times $T^\nu$ such that exact solutions $u^\nu$ to the Navier-Stokes equations satisfy 
$$
\| u^\nu(0) - U_s(0) \|_{H^s} \le \nu^N,
$$
but 
$$\| u^\nu (T^\nu) - U_s(T^\nu)\|_{L^\infty} \ge \sigma_0,$$
for time sequences $T^\nu \to 0$, as $\nu \to 0$. 

\item There is a source $f^\nu$ that is sufficiently small in $L^1(\RR_+;L^\infty(\Omega))$ and is exponentially localized within the boundary layer of size $\nu^{3/4}$ so that the following holds: there is a positive constant $\sigma_0$ so that the unique solution $u^\nu$ of the Navier-Stokes equations, with source $f^\nu$ in the momentum equation and with the initial data $u^\nu_{\vert_{t=0}} = {v_S^1}_{\vert_{t=0}}$, must satisfy 
$$ \| u^\nu(T^\nu) - v_S^1(T^\nu) \|_{L^\infty} \ge \sigma_0,$$ 
for time sequences $T^\nu \to 0$, as $\nu \to 0$.

 \end{itemize}

\end{theorem}

A more precise result is given in the end of the paper, where initial perturbations and sources for the stability of the boundary sublayer are more explicit.

%%%%%%%%%%%%%%%%%%%%%%%%%%

\section{Construction of an approximate solution}

%%%%%%%%%%%%%%%%%%%%%%%%%%

%%%%%%%%%%%%%%%%%%%%%%%

\subsection{A first scaling}

%%%%%%%%%%%%%%%%%%%%%%%

We first rescale time and space according to the classical change of variables
$$
T = {t \over \sqrt{\nu}}, \quad
X = {x \over \sqrt{\nu}}, 
\quad
Y = {y \over \sqrt{\nu}}.
$$
The Navier-Stokes equations \eqref{NS1}-\eqref{NS2} are invariant under this scaling, except the viscosity coefficient which is now $\sqrt{\nu}$ instead of $\nu$. For the rest of the paper, we shall work with the scaled Navier-Stokes equations with above scaled variables. For sake of presentation, we write $t,x,y$ in place of  $T$, $X$ and $Y$, respectively. 

Let $U(y)$ be the inviscid unstable shear flow, and let $U_s(\sqrt \nu t,y)$ be the corresponding time-dependent shear flow. Our goal is to construct an approximate solution to the Navier-Stokes equations that exhibits instability. Let us introduce 
$$ v = u - U_s$$
in which $u$ is the genuine solution to the Navier-Stokes equations. Then, $v$ solves 
\beq \label{NSI1}
\partial_t v + (U\cdot \nabla) v + (v \cdot \nabla) U + (v \cdot \nabla) v + \nabla p =  \sqrt \nu\Delta v + \sqrt \nu S v
\eeq
\beq \label{NSI2}
\nabla \cdot v = 0,
\eeq
in which the linear operator $Sv$ is defined by 
\begin{equation}\label{def-Sv}
 Sv := \nu^{-1/2} [U_s (\sqrt \nu t)  - U] \cdot \nabla v + \nu^{-1/2} v \cdot \nabla [U_s (\sqrt \nu t)  - U] .
\end{equation}

We shall establish the instability
in four steps. First, we construct an approximate solution that exhibits the instability, starting from the maximal linear growing mode of Euler equations. 
We then construct an approximate viscous solution to the Navier-Stokes equations in the large scale of size $\sqrt \nu$, without correcting the no-slip boundary condition. That is, the approximate solution satisfies only the zero normal velocity condition 
\beq \label{NSI3}
u \cdot n = 0 \quad \hbox{on} \quad \partial \Omega .
\eeq
In the third step, we approximately correct the nonzero boundary value caused by the previous step. This leads to an instability, which is a 
solution of Navier-Stokes equations, except for a small error term which is localized in a layer of thickness of order $\nu^{3/4}$.
The remainder of the paper is devoted to the study of the stability of this approximate solution.

%%%%%%%%%%%%%%%%%%%%%%%

\subsection{Linear instability}

%%%%%%%%%%%%%%%%%%%%%%%

By assumption, $U(y)$ is spectrally unstable for Euler equations. That is, there exists a solution of linearized Euler equations, namely equations (\ref{NSI1})-(\ref{NSI3})
with $\nu = 0$, of the form
\begin{equation}\label{gr-soln}
u_e^0 = \nu^N \Re (u_e e^{\lambda t} ) 
\end{equation}
where $u_e$ is a smooth, divergence free, vector field and $\Re \lambda > 0$. Since the unstable spectrum of the linearized Euler equations around a shear flow consists of only unstable eigenvalues, we assume that $\lambda$ is the maximal unstable eigenvalue. 
 Furthermore, as $u_e$ is divergence free, it can be written under the form
$$
u_e = \nabla^\perp ( \psi_e e^{i \alpha x}) .
$$
Here, $\alpha$ denotes the wave number of the Fourier transform and $\psi_e$ is the corresponding stream function, both solving the corresponding Rayleigh equation:
\beq \label{Ray}
(U - c) (\partial_y^2 - \alpha^2) \psi_e = U'' \psi_e
\eeq
with boundary conditions $\psi_e(0) = \lim_{y \to + \infty} \psi_e = 0$, with $c = -\lambda / i\alpha$. Since $\psi_e$ is a solution of an elliptic equation, it is real analytic. As a consequence, the unstable eigenfunction $u_e$ is entire in $x$ and holomorphic on $y$.

In addition, it follows that the $L^p$ norm of $u_e^0$ behave like $\nu^N e^{\Re \lambda t}$. Precisely, there are positive constants $c_0,c_1$ so that 
\begin{equation}\label{Egr} c_0 \nu^N e^{\Re \lambda t} \le \| u^0_e(t) \|_{L^p} \le c_1 \nu^N e^{\Re \lambda t}, \qquad 1\le p \le \infty \end{equation} 
in which $\|\cdot \|_{L^p}$ denotes the usual $L^p$ norm. 
Let us introduce the instability time  $T^\star$, defined by
\begin{equation}\label{def-T}
T^\star = - N {\log \nu \over \Re \lambda} 
\end{equation}
and time $T^\star_\theta$, for any $\theta \ge 0$, defined by
\begin{equation}\label{def-Ttheta}
T^\star_\theta = - (N - \theta) {\log \nu \over \Re \lambda}. 
\end{equation}
We observe that by \eqref{Egr}, $\| u_e^0(T^\star_\theta) \|_{L^p}$ is exactly of order
$\nu^\theta$. In order to get order one instabilities for Prandtl's layers, it is necessary to construct a solution of Navier
Stokes equations up to the time $T^\star$, or at least to $T^\star - \tau$ for some possibly large, but fixed $\tau$, in the inviscid limit.
However, this appears to be very difficult due to the appearance of a viscous boundary sublayer of order $\nu^{1/4}$. The presence of such a sublayer causes the (viscous) approximate solution to have a large gradient of order $\nu^{N-\frac14} e^{\Re \lambda t}$. For this reason, the approach introduced in \cite{Grenier} stops at the time $T_{1/4}^\star$. After this time, energy estimates cannot be fulfilled. The aim of this construction is to investigate what appends between $T^\star_{1/4}$ and $T^\star_0$.

%%%%%%%%%%%%%%%%%%%%%%%

\subsection{Construction of an "inviscid" nonlinear instability}

%%%%%%%%%%%%%%%%%%%%%%%

In this section, we build an approximate solution of \eqref{NSI1}-\eqref{NSI3}, starting from $u_e^0$. We stress that this solution only satisfies the boundary condition \eqref{NSI3} for Euler solutions.
Precisely, we construct solutions of the form 
\begin{equation}\label{def-uapp-e} u_e^\app = \nu^N \sum_{j=0}^M \nu^{j/2} u_e^j.\end{equation}
For sake of simplicity, we take $N$ to be a (sufficiently large) integer. Plugging this Ansatz into \eqref{NSI1} and matching order in $\nu$, we are led to solve 

\begin{itemize}

\item for $j =0$: $u^0_e$ is the growing solution defined in \eqref{gr-soln}. 
 
\item for $0<j\le M$: 
\begin{equation}\label{eqs-omegaj} 
\begin{aligned}
\partial_t u_e^j + (U \cdot \nabla) u_e^j + (u_e^j \cdot \nabla) U + \nabla p  &= R_j, 
\\
\nabla \cdot u_e^j &=0,
\\
u_e^j \cdot n &=0, \qquad \mbox{on}\quad \partial \Omega,
\end{aligned}
\end{equation}
\end{itemize}
together with zero initial data. Here, the remainders $R_j$ is defined by 
$$ 
R_j = Su_e^{j-1} + \Delta u_e^{j-1}+ \sum_{k + \ell + 2N = j} u_e^k \cdot \nabla u_e^\ell.
$$
As a consequence, $u_e^\app$ appropriately solves the Navier-Stokes equations \eqref{NSI1}-\eqref{NSI2}, with the Euler's boundary condition \eqref{NSI3}, leaving an error of the approximation $E_e^\app$, defined by 
$$
\begin{aligned}
E_e^\app
&= \nu^{N+\frac{M+1}{2}} ( S u_e^M + \Delta u_e^M)+  \sum_{k+ \ell> M+1 -2 N; 1\le k,\ell \le M} \nu^{2N+ \frac{k+\ell}{2}} u_e^k \cdot \nabla u_e^\ell .
\end{aligned}
$$
Note that at each step, $u_e^k$ is a solution of linearized Euler equations around $U$, 
with a source term consisting of solutions constructed in the previous steps. Since the linearized Euler problem is well-posed, 
$u_e^j$ is uniquely defined. In addition, by letting $L$ be the linearized Euler operator around $U$, 
there holds the uniform semigroup estimate 
$$
 \| e^{Lt} u_e \|_{H^s} \le C_\beta e^{(\Re \lambda + \beta) t} \|u_e \|_{H^{s+2}} , \qquad \forall~ \beta >0, \qquad \forall ~t\ge 0,
 $$
for all $s\ge 0$. The loss of derivatives in the above semigroup can be avoided by studying the resolvent solutions to the Euler equations 
or the Rayleigh equations (similarly, but much simpler, to what is done for linearized Navier-Stokes; see, for instance, \cite{GrN1,GrN2}). 

By induction, by using the semigroup estimates, it is then straightforward (e.g., \cite{Grenier}) to prove that 
\begin{equation}\label{induction-j} \| u_e^j\|_{H^{s+4M-4j}} \le C_{j,s} e^{(1+ \frac{j}{2N}) \Re \lambda t} \end{equation}
for any $s\ge 0$. 
As a consequence, the function $u_e^\app$ defined as in \eqref{def-uapp-e} approximately solves the Navier-Stokes equations in the following sense: 
\begin{equation}\label{eqs-omegaj} 
\begin{aligned}
\partial_t u_e^\app + (U_s + u_e^\app) \cdot \nabla u_e^\app + u_e^\app \cdot \nabla U_s + \nabla p  &= \sqrt \nu \Delta u_e^\app + E_e^\app, 
\\
\nabla \cdot u_e^\app &=0,
\\
u_e^\app \cdot n &=0, \qquad \mbox{on}\quad \partial \Omega,
\end{aligned}
\end{equation}
in which $U_s = U_s(\sqrt \nu t,y)$. In addition, as long as $\nu^N e^{\Re \lambda t}$ remains bounded, there hold 
\begin{equation}\label{bd-uapp-e}
\| u_e^{\app} \|_{H^s} \le C \nu^{N} e^{\Re \lambda t}, \qquad \| E_e^\app \|_{H^s} \le  C \Big(  \nu^{N} e^{\Re \lambda t} \Big)^{1+ \frac{M+1}{2N}}. 
\end{equation}
Again, we stress that the approximate solution $u_e^\app$ does not satisfy the no-slip boundary condition on $\partial \Omega$, but the condition \eqref{NSI3} on the normal component of velocity. 
Note that, in particular, for any $\theta > 0$ and for $t \le T^\star_\theta$, there holds
$$
\| E^\app_e \|_{L^2} \le C \nu^{\theta P},
$$
which can be made arbitrarily small if $P = 1+ \frac{M+1}{2N}$ is chosen large enough.
Roughly speaking, $u_e^\app$ describes the "large scale" instability, which we shall introduce in the next section.

\begin{remark} The approximate solution $u_e^\app$ is in fact holomorphic on $\Omega$. Indeed, it suffices to prove the claim that $u_e^j$ is a linear combination of functions of the form $\nabla^\perp (\psi e^{i \beta x} )$. Indeed, by construction, the claim holds for $u_e^0$. Assume that the claim holds for $j\ge 0$. Then, in particular, the source $R_j$ is holomorphic and also a linear combination of functions of the form $\nabla^\perp (\psi e^{i \beta x} )$. Taking Fourier-Laplace transform, we get that
$u_e^{j+1}$ is a sum of solutions of Rayleigh equations, and is therefore holomorphic.
\end{remark}

%%%%%%%%%%%%%%%%%%

\subsection{Large scale behavior}

%%%%%%%%%%%%%%%%%%

We now look for a corrector $\widetilde  u_e$ of $u_e^{\app}$ which kills the large scale error term $E^\app_e$. 
Precisely, we construct the corrector $\widetilde  u_e$ so that 
\begin{equation}\label{def-uL}
u_L: = U_s(\sqrt \nu t,y) + u_e^\app + \widetilde u_e
\end{equation}
is an exact solution to the Navier-Stokes equations, without taking care of the no-slip boundary condition. Indeed, we shall replace the no-slip condition by a Navier boundary condition, which allows us to derive uniform bounds on vorticity. The no-slip boundary condition will then be recovered in the next section. 

The large scale corrector $\widetilde u_e$, defined as in \eqref{def-uL}, solves
$$
\partial_t \widetilde  u_e + u_L \cdot \nabla \widetilde  u_e + \widetilde  u_e \cdot \nabla (U_s  + u_e^{\app}) + \nabla p - \sqrt \nu \Delta \widetilde  u_e = - E_e^\app,
$$
$$
\nabla \cdot \widetilde  u_e = 0,
$$
with zero initial data $\widetilde u_e=0$ at $t=0$, and with the following Navier boundary conditions
$$ \widetilde u_e \cdot n =0, \qquad (D\widetilde u_e )n \cdot \tau =0, $$
on $\partial \Omega$. Here, $Du = \frac12 (\nabla u + (\nabla u)^{tr})$. On the flat boundary, the above Navier boundary conditions in particular yield $\widetilde \omega_e =0$ on $\partial\Omega$. 
We stress that $u_L$ does not satisfy the no-slip boundary condition. However, it describes the large scale behavior of the main Prandtl's boundary layer.

By energy estimates, using the fact that $\| \nabla u_e^{\app}\|_{L^\infty}$ is bounded by $\nu^N e^{\Re \lambda t}$, and using the zero boundary condition on the normal component of velocities, we get
$$
\frac12 \frac{d}{dt} \| \widetilde  u_e \|_{L^2}^2 \le C  (1 + \nu^N e^{\Re \lambda t}) \| \widetilde  u_e \|_{L^2}^2 + C  \Big(  \nu^{N} e^{\Re \lambda t} \Big)^{2+ \frac{2(M+1)}{2N}} .
$$
Hence, as long as $\nu^N e^{\Re \lambda t}$ remains bounded (or equivalently, $t \le T^\star$),
this yields 
\begin{equation}\label{bd-tue}
\| \widetilde  u_e (t) \|_{L^2} \le C  \Big(  \nu^{N} e^{\Re \lambda t} \Big)^{1+ \frac{M+1}{2N}} .
\end{equation}
In particular, for any $\theta > 0$ and for $t \le T^\star_\theta$, there holds
$$
\| \widetilde  u_e \|_{L^2} \le C \nu^{\theta P}
$$
which can be arbitrarily small with respect to $\nu$, provided $P\gg1$. Similarly, since $x$-derivatives of $\widetilde u_e$ satisfy the same type of boundary conditions, there hold 
$$
\| \partial_x^k \widetilde  u_e \|_{L^2} \le C  \Big(  \nu^{N} e^{\Re \lambda t} \Big)^{P} , \qquad \forall k\ge 0.  
$$

In addition, the standard elliptic estimates on $\TT\times \RR_+$ yield $\| \widetilde u_e \|_{L^\infty} \lesssim \| \widetilde \omega_e\|_{L^\infty}$. To bound the vorticity $\widetilde \omega_e$, we write 
$$
\partial_t \widetilde  \omega_e + u_L \cdot \nabla \widetilde  \omega_e  - \sqrt \nu \Delta \widetilde  \omega_e = - \widetilde  u_e \cdot \nabla (\omega_s  + \omega_e^{\app})- \nabla \times E_e^\app
$$
with $\widetilde \omega_e =0$ on the boundary. The Maximum Principle for the transport-diffusion equation, together with \eqref{bd-uapp-e}, yields 
$$
\begin{aligned}
 \| \widetilde \omega_e (t)\|_{L^\infty } 
 &\le \int_0^t \Big[ \| \widetilde  u_e \cdot \nabla \omega_s\|_{L^\infty}  +\| \widetilde  u_e \cdot \nabla \omega_e^{\app}\|_{L^\infty} +  \|\nabla \times E_e^\app\|_{L^\infty}\Big] \; ds 
\\
 &\le C\int_0^t \Big[ \|\widetilde  u_{e,2}\partial_y\omega_s\|_{L^\infty}  +\nu^N e^{\Re \lambda s}\| \widetilde  \omega_e(s)\|_{L^\infty} +  \Big(  \nu^{N} e^{\Re \lambda s} \Big)^{P}
 \Big] \; ds .
  \end{aligned}$$
Writing $\widetilde  u_{e,2} = \int_0^y \partial_y\widetilde  u_{e,2} \; dy$, we have $$\begin{aligned}
 \|\widetilde  u_{e,2}\partial_y\omega_s\|_{L^\infty} \le  \| y \partial_y \omega_s\|_{L^\infty} \| \partial_x \widetilde u_{e,1}\|_{H^1_xL_y^2} \le C 
 \Big(  \nu^{N} e^{\Re \lambda t} \Big)^{P}.
 \end{aligned}$$ 
Thus, as long as $\nu^N e^{\Re \lambda t}$ remains sufficiently small (or equivalently, $t \le T^\star-\tau$ for large $\tau$), we obtain at once 
\begin{equation}\label{bd-twe}
\begin{aligned}
 \| \widetilde \omega_e (t)\|_{L^\infty } 
 &\le C\Big(  \nu^{N} e^{\Re \lambda t} \Big)^{P}.
  \end{aligned}\end{equation}
This yields the same bound for velocity $\widetilde u_e$ and $\partial_x^k \widetilde u_e$ in $L^\infty$, for $k\ge 0$.

%%%%%%%%%%%%%%%%%%%%%

\subsection{Sublayer correction}\label{sec-sublayer}

%%%%%%%%%%%%%%%%%%%%%

It remains to correct the no-slip boundary condition of the (exact) solution $u_L$. To this end, we introduce the sublayer correction $u_S$, solving the Navier-Stokes equation  
\begin{equation}\label{eqs-uS}
\begin{aligned}
\partial_t u_S + (u_L \cdot \nabla) u_S + (u_S \cdot \nabla) u_L + (u_S \cdot \nabla) u_S + \nabla p  &=  \sqrt \nu \Delta u_S,
\\
\nabla \cdot u_S &= 0,
\end{aligned}\end{equation}
together with the inhomogenous boundary condition 
\begin{equation}\label{bdry-uS}u_S = - u_L = - u_e^\app - \widetilde u_e, \qquad  \mbox{on}\quad \partial \Omega,\end{equation}
in which $u_e^\app+\widetilde u_e$ is of order $\nu^N e^{\Re \lambda t}$; see \eqref{bd-uapp-e} and \eqref{bd-twe}. Observe that 
$$u = u_S + u_L$$ 
is an exact solution of the genuine Navier-Stokes equations \eqref{NS1}-\eqref{NS2}, with the no-slip boundary condition \eqref{NS3}.
As we will see, $u_S$ describes the "small structures" of $u$, namely its viscous boundary sublayer.

To solve \eqref{eqs-uS}, let us first consider the following simplified equations 
$$
\partial_t u_S^1 + (U_s \cdot \nabla) u_S^1 + (u_S^1 \cdot \nabla) U_s + \nabla p =\sqrt \nu \Delta u_S^1,
$$
$$
\nabla \cdot u_S^1 = 0,
$$
with the boundary condition $
u_S^1 = - u_L
$
on $\partial \Omega$. Note that $u_S^1$ has a boundary layer behavior, with a small scale of order $\nu^{1/4}$ in $y$.
As a consequence, as $U_s$ is tangential to the boundary and is of order $O(y)$ for small $y$, the convection terms $(U_s \cdot \nabla) u_S^1$
and $(u_S^1 \cdot \nabla) U_s$ are of order $O(\nu^{1/4})$ smaller than $u_S^1$. Thus, the convection terms might be moved into the next order and $u_S^1$ may be
approximated by $v_S^1$, a solution of the linear Stokes equation
\begin{equation}\label{Stokes}
\partial_t v_S^1 + \nabla p = \sqrt \nu \Delta v_S^1 ,
\qquad \nabla \cdot v_S^1 = 0,
\end{equation}
with the same boundary condition $v_S^1 = - u_L$. 

The Stokes problem can be solved explicitly by introducing the stream function
$\phi_c$ defined through
$$
v_S^1 := - \nabla^\perp \phi_c .
$$
Starting with $v_S^1$, we can construct an approximation of $u_S$. Again, our construction is inductive. For $k \ge 2,$ we iteratively construct $v_S^k$, solving the following Stokes problem 
$$
\partial_t v_S^k + \nabla p - \sqrt \nu \Delta v_S^k = - Q_k,
$$
$$
\nabla \cdot v_S^k = 0
$$
with the zero Dirichlet boundary condition on $v_S^k$ and zero initial data. Here, the remainder $Q_k$ is defined by 
$$
Q_k = (u_L \cdot  \nabla) v_S^{k-1} + (v_S^{k-1} \cdot \nabla ) u_L  + \sum_{j+\ell = k} (v_S^j \cdot  \nabla) v_S^{\ell} .
$$
We then set 
$$
u_S^{\app} = \sum_{k =1 }^M v_S^k
$$
where $M$ is arbitrarily large. By construction
$u_S^\app$  approximately solves the Navier-Stokes equations \eqref{eqs-uS}, leaving an error $R^\app_S$ in the momentum equation. It is then straightforward to prove that 
\begin{equation}\label{bd-uSapp} | u_S^\app(t,x,y)| + |\partial_x u_S^\app(t,x,y)|\le C \nu^N e^{\Re \lambda t} e^{-\beta y /\nu^{1/4}}\end{equation}
for some positive constant $\beta$, and the remainder $R_S^\app$ satisfies 
\begin{equation}\label{bd-RSapp}|R_S^\app(t,x,y)| \le C e^{-\beta y /\nu^{1/4}} \Big( \nu^N e^{\Re \lambda t}\Big)^P ,\end{equation}
for some positive $P$, which can be taken to be arbitrarily large (for large enough $M$ in the construction of the approximate solution).

%%%%%%%%%%%

\subsection{Approximate solution}

%%%%%%%%%%%

We are ready to conclude the construction of an approximate solution. Indeed, introduce 
\begin{equation}\label{def-uapp}
u^{\app} = u_L + u_S^{\app}  = U_s(\sqrt \nu t,y) + u_e^\app + \widetilde u_e + u_S^\app
\end{equation}
with  $u_e^\app, \widetilde u_e,$ and $u_S^\app$ constructed in the previous subsections. Then, $u^\app$ approximately solves the nonlinear Navier-Stokes equations in the following sense
\begin{equation}\label{NS-app} 
\begin{aligned}
\partial_t u^\app + u^\app \cdot \nabla u^\app + \nabla p  &= \sqrt \nu \Delta u^\app + R_S^\app, 
\\
\nabla \cdot u^\app &=0,
\\
u^\app  &=0, \qquad \mbox{on}\quad \partial \Omega,
\end{aligned}
\end{equation}
with the remainder $R^\app_S$ satisfying \eqref{bd-RSapp}. We stress that the remainder $R_S^\app$ is exponentially localized near the boundary with thickness of order $\nu^{1/4}$. 
 
Let us detail the structure of this approximate solution. By construction, we recall that 
$$\| \partial_x^k u_e^\app(t)\|_{L^\infty} \le C\nu^N e^{\Re \lambda t}, \qquad  \|\partial_x^k \widetilde u_e\|_{L^\infty} \le C \Big( \nu^N e^{\Re \lambda t} \Big)^P,$$
for $k\ge 0$, and recall the estimate \eqref{bd-uSapp} for $u_S^\app$. This in particular yields 
$$ \| \partial_x u^\app \|_{L^\infty} \le C \nu^N e^{\Re \lambda t} .$$
Using divergence-free condition, we thus get the same bound for $\partial_y u^\app_2$, and hence we obtain the following pointwise bound
\begin{equation}\label{uapp-y}
| u^\app_2(t,x,y) | \le C \nu^N e^{\Re \lambda t} y .
\end{equation}

Let us give a lower bound on the approximate solution. By view of \eqref{Egr}, there exists some positive constant $c_2$ 
so that \begin{equation}\label{lower-bd}
\| u^\app - U_s(\sqrt \nu t,\cdot)\|_{L^\infty} \ge c_2 \nu^N e^{\Re \lambda t},
\end{equation}
for all $t\ge 0$, as long as $\nu^N e^{\Re \lambda t}$ remains sufficiently small (independent of $\nu$).

%

%%%%%%%%%%%%%%%%%%%%%%%%%

\section{Sublayer behavior}

%%%%%%%%%%%%%%%%%%%%%%%%%

%%%%%%%%%%%

\subsection{Link between sublayer and Prandtl layer}

%%%%%%%%%%%

Let $u^\nu$ be the genuine solution to the Navier-Stokes equations, and let $u^\app$ be the approximate solution constructed in the previous section. Set
$$
v = u^\nu - u^\app .
$$
It follows that $v$ solves 
\beq \label{NS1v}
\begin{aligned}
\partial_t v + (u^\app +v)\cdot \nabla v + v \cdot \nabla u^\app + \nabla p &= \sqrt\nu\Delta v  + R_S^\app,
\\
\nabla \cdot v &= 0,
\end{aligned}\eeq
with the zero Dirichlet boundary condition and with an exponentially localized remainder $R_S^\app$ satisfying 
\begin{equation}\label{bd-RSapp1}|R_S^\app(t,x,y)| \le C e^{-\beta y /\nu^{1/4}} \Big( \nu^N e^{\Re \lambda t}\Big)^P .\end{equation}

As the source term $R_S^\app$ is located in the sublayer, we expect that $v$ is also located in the sublayer, provided
the vertical transport remains small. 
Note that (\ref{NS1v}) describes the behavior of a boundary sublayer of size $\nu^{1/4}$ (and hence of size $\nu^{3/4}$ in the original variables). Let us make yet another change of variables:
$$
X = {x \over \nu^{1/4}}, \quad Y = {y \over \nu^{1/4}}, \quad T = {t \over \nu^{1/4}} .
$$

Then, in these new variables, (\ref{NS1v}) becomes the (same) Navier-Stokes equations, with viscosity $\nu^{1/4}$, near an approximate
solution $u^\app$ which exhibits a boundary layer behavior. Precisely, $u^\app$ is of the form
$$u^\app =  u_L(\nu^{1/4}T, \nu^{1/4}X, \nu^{1/4} Y) + u_S^\app(\nu^{1/4} T, \nu^{1/4} X,Y)$$
in which the leading term $v_S^1$ in $u_S^\app$ solves the Stokes problem \eqref{Stokes}. The boundary type approximate solution $u^\app$ is very close to a Prandtl's profile, except that
\begin{itemize}

\item It is slowly evolving in the $X$ direction, on sizes of order $\nu^{1/4}$

%\item It is slowly evolving in the $Y$ direction, on sizes of order $\nu^{1/4}$

\item There is a small upward velocity, of order $\nu^{1/4}$.

\end{itemize}

Roughly speaking, equation  (\ref{NS1v}) describes the stability of approximate boundary layer solutions,
which are small amplitude and slow modulations of pure shear layers.
It is very likely that if Prandtl layers are stable, so is  (\ref{NS1v}), since it is reasonable to believe that
any proof of stability for Prandtl layers would bear small perturbations and slow spatial modulations. As we will see in the next
paragraph, this belief appears to be false.

%%%%%%%%%%

\subsection{Stability of the sublayer}

%%%%%%%%%%

Let us assume that the sublayer is nonlinearly stable in $L^\infty$; namely, we assume either $\|v(t)\|_{L^\infty}$ remains sufficiently small or 
$$
\| v (t) \|_{L^\infty} \le C \Big( \nu^{N} e^{\Re \lambda t} \Big)^{1+\beta}
$$
for $t\le T^\star $ and for some $\beta>0$. Note that this notion of stability is very weak, since we  expect $\beta \ge 1$.
Then, the Prandtl layer is nonlinearly unstable, since at $t = T^\star - \tau$ for $\tau$ large enough, $\nu^N e^{\Re \lambda t}$ remains sufficiently small and hence \eqref{lower-bd} yields
$$
\| u^\nu - U_s \|_{L^\infty} \ge \| u^\app - U_s \|_{L^\infty}  - \| v \|_{L^\infty} \ge \sigma_0 > 0
$$
for some positive (and small) constant $\sigma_0$ (independent on $\nu$). The main theorem is proved.

%%%%%%%%%%%%%%%%%%%%%%%%%%%%%%%%%%%%%%%%%%%%

\end{document}